\documentclass[12pt,a4paper]{article}

\usepackage{amssymb,amsfonts}       % AMS Fonts Family
\usepackage{epsfig}
\usepackage{amsmath}
\usepackage{amsthm}
\usepackage{latexsym}

\newcommand{\vol}{\mathop{\rm Vol}\nolimits}

\newcommand{\pr}{\mathop{\rm pr}\nolimits}

\numberwithin{equation}{section}
\swapnumbers

\def\?{\marginpar{$\bullet\bullet\bullet$}}

\title{\bf Gravity model for topological features \\ on a
cylindrical manifold}
\author{Igor Bayak}
\date{}

\begin{document}
\maketitle
\begin{abstract}

A model aimed at understanding quantum gravity in terms of
Birkhoff's approach is discussed. The geometry of this model is
constructed by using a winding map of Minkowski space into a
$\mathbb{R}^{3} \times S^{1}$-cylinder. The basic field of this
model is a field of unit vectors defined through the velocity
field of a flow wrapping the cylinder. The degeneration of some
parts of the flow into circles (topological features) results in
inhomogeneities and gives rise to a scalar field, analogous to the
gravitational field. The geometry and dynamics of this field are
briefly discussed. We treat the intersections between the
topological features and the observer's 3-space as matter
particles and argue that these entities are likely to possess some
quantum properties.
\end{abstract}

\section{Introduction}

In this paper we shall discuss a mathematical construction aimed
at understanding quantum gravity in terms of Birkhoff's twist
Hamiltonian diffeomorphism of a cylinder \cite{BIR}. We shall also
use the idea of compactification of extra dimensions due to Klein
\cite{KLE}. To outline the main idea behind this model in a very
simple way, we can reduce the dimensionality and consider
the dynamics of a vector field defined on a 2-cylinder
$\mathbb{R}^{1} \times S^{1}$. For this purpose we can use
the velocity field $u(x,\tau)$ of a two-dimensional flow of ideal
incompressible fluid moving through this manifold.

Indeed, the dynamics of the vector field $u(x,\tau)$ with the
initial condition $u(x,0)$ is defined by the evolution
equation
\begin{equation}\label{EvoldynCyl}
  \delta\int_{\Delta \tau}\int_{\Delta x}dx\wedge
  u(x,\tau)d\tau\rightarrow 0,
\end{equation}
where we use the restriction of the vector field onto an arbitrary
cylinder's element; $\Delta \tau$ is
the evolution (time) interval, and $\Delta x$ is an arbitrary
segment of the cylinder's element. In other words, we assume the
variation of the integral of the mass carried by the flow through
the segment during a finite time interval to be vanishing. That
is, as a result of the field evolution, $u(x,0)\rightarrow
u(x,\infty)$, the functional of the flow mass approaches to its
maximal value. If, at the initial moment of time, the regular
vector field $u(x,0)$ corresponds to a unit vector forming an
angle $\varphi$ with the cylinder's element, then the evolution of
this field is described by the equation
\begin{equation}\label{EvoldynReg}
  \delta\int_{\Delta \tau}\int_{\Delta x}dx\wedge
  u(x,\tau)d\tau= \delta\int_{\Delta \tau}\int_{\Delta x}\sin\varphi(\tau)
  dxd\tau=\cos\varphi(\tau)\Delta \tau\Delta x\rightarrow 0.
\end{equation}
Therefore, the case of $\varphi(0)=0$ corresponds to
the absolute instability of the vector field. During its evolution,
 $u(x,0)\rightarrow u(x,\infty)$, the field is relatively stable at
 $0<\varphi(\tau)<\pi/2$, achieving the absolute stability
  at the end of this evolution, when $\varphi(\infty)=\pi/2$.
If, additionally, we fix the vector field $u(x,\tau)$
at the endpoints of the segment $\Delta x$ by imposing
some boundary conditions on the evolution equation
\ref{EvoldynCyl}, we would get the following dynamical
equation:
\begin{equation}\label{DynCyl}
  \delta\int_{\Delta \tau}\int_{\Delta x}dx\wedge
  u(x,t)d\tau = 0.
\end{equation}
Let some flow lines of the vector field $u(x,\tau)$ be
degenerated into circles (topological features) as a result
of the absolute instability of the field and
fluctuations during the initial phase of its evolution.
Since the dynamics of such topological features is described
by \ref{DynCyl}, the features would tend to move towards that
side of $\Delta x$ where the field $u(x,\tau)$
is more stable. Thus, the topological features serve as
attraction points for each other and can be used for modelling
matter particles (mass points).

We must emphasise that the plane $(x,\tau)$, in which our
variational equations are defined, has the Euclidean metric. That
is, in the case of the Euclidean plane $(x,\phi)$ wrapping over a
cylinder we can identify the azimuthal parameter $\phi$ with the
evolution parameter $\tau$. By choosing the observer's worldline
coinciding with a cylinder's element we can speak of a classical
limit, whereas by generalising and involving also the azimuthal
(angular) parameter we can speak of the quantisation of our model.
So, when the observer's
worldline is an arbitrary helix on the cylinder, the
variational equation \ref{DynCyl} reads
\begin{equation}\label{MinDynCyl}
  \delta\int_{\Delta x_{0}}\int_{\Delta x_{1}}dx_{1}\wedge
  g(x)dx_{0} = 0,
\end{equation}
where the varied is the vector field $g(x)$ defined on the
pseudo-Euclidean plane $(x_{0},x_{1})$ oriented
in such a way that one of its isotropic lines covers the
cylinder-defining circle and the other corresponds to
a cylinder's element. In this case we can speak of
a relativistic consideration. If the observer's worldline
corresponds to a curved line orthogonal to the flow lines
of the vector field  $g(x)$, where $g^{2}(x)>0$,
then we have to use the variational equation defined on
a two-dimensional pseudo-Riemann manifold  $M$ induced
by the vector field $g(x)$, namely,
\begin{equation}\label{RimDynCyl}
  \delta\int_{\Delta M}g^{2}(x')\sqrt{-\det g_{ij}}\, dx'_{0}\wedge
  dx'_{1} = 0,
\end{equation}
where $\Delta M=\Delta x'_{0}\times\Delta x'_{1}$  is an arbitrary
region of the manifold $M$; $x'_{0}(\phi)$ is the flow line of the
vector field  $g(x)$ parameterised by the angular coordinate
$\phi$; $x'_{1}(r)$ is the spatial coordinate on the cylinder
(orthogonal to the observer worldline) parameterised by the
Euclidean length $r$; $g_{ij}$ is the Gram matrix corresponding to
the pair of tangent vectors
$\left(\frac{dx'_{0}}{d\phi},\frac{dx'_{1}}{dr}\right)$. In this
case the dynamics of the vector field is described through the
geometry of its flow lines \cite{Johnson,Aminov,Yampolsky}.

Thus, we can say that our approach to the dynamics of the vector
field is based on maximisation of the mass carried by the flow
\cite{gibbons02,schiller05}, which is not exactly what is
typically used in the ergodic theory
\cite{birkhoff31,hopf37,sinai76}. However, this principle is
likely to be related to the the minimum principle for the velocity
field \cite{reiser96,bejan04,montassar05},  which is a special
case of the more general principle of minimum or maximum entropy
production \cite{jaynes57,hamann69}.

Before a more detailed discussion of this model we have to make a
few preliminary notes. First, throughout this paper we shall use a
somewhat unconventional spherical coordinates. Namely, latitude
will be measured modulo $2\pi$ and longitude -- modulo $\pi$. In
other words,  we shall use the following spherical ($\rho$,
$\varphi$, $\theta_{1},\ldots,\theta_{n-2}$) to Cartesian
($x_{1},\ldots,x_{n}$) coordinate transformation in
$\mathbb{R}^{n}$:
\begin{align*}
&x_{1}=\rho\cos\varphi, \\
&x_{2}=\rho\sin\varphi\cos\theta_{1}, \\
&x_{3}=\rho\sin\varphi\sin\theta_{1}, \\
&\hspace{0.1cm} \dots \dots \dots \dots \dots \dots \dots \\
&x_{n-1}=\rho\sin\varphi \hspace{0.1cm} \ldots
\hspace{0.1cm} \sin\theta_{n-3}\cos\theta_{n-2}, \\
&x_{n}=\rho\sin\varphi \hspace{0.1cm} \ldots
\hspace{0.1cm} \sin\theta_{n-3}\sin\theta_{n-2},
\end{align*}
where $0\leq\rho<\infty$, $0\leq\varphi<2\pi$ and $0\leq\theta_{i}<\pi$.
We shall also be interpreting the projective space $RP^{n}$ as the space
of centrally symmetric lines in $\mathbb{R}^{n+1}$, that is, as a quotient space
$\mathbb{R}^{n+1}\backslash \{0\}$ under
the equivalence relation $x\sim rx$, where
$r\in \mathbb{R}\backslash \{0\}$.

\section{The geometry of the model}

We can describe the geometry of our model in terms of
the mapping of the Euclidean
plane into a 2-sphere, $S^{2}$, by winding the former around
the latter. We can also use similar winding maps for the
pseudo-Euclidean plane into a cylinder, $\mathbb{R}\times S^{1}$,
or a torus, $S^{1}\times S^{1}$. More formally this could be
expressed in the following way \cite{bialy92}.
Take the polar coordinates $(\varphi,\rho)$ defined on the
Euclidean plane and the spherical coordinates $(\theta,\phi)$
on a sphere. We can map the
Euclidean plane into sphere by using the congruence
classes modulo $\pi$ and $2\pi$. That is,
\begin{equation}\label{sphere}
  \theta=|\varphi|\mod \pi,\quad \phi=|\pm\pi\rho|\mod 2\pi,
\end{equation}
where the positive sign corresponds to the interval
$0\leq\varphi<\pi$ and negative -- to the interval
$\pi\leq\varphi<2\pi$.  If the projective
lines are chosen to be centrally symmetric then the Euclidean
plane can be generated as the product $RP^{1}\times \mathbb{R}$.
Here the components of $\mathbb{R}$ are assumed to be Euclidean,
i.e., rigid and with no mirror-reflection operation allowed.
Similarly,  we can define a space based on unoriented lines in the
tangent plane to the sphere. Therefore, the sphere can be
generated by the product $RP^{1}\times S^{1}$, the opposite points
of the circle being identified with each other. In this
representation all centrally symmetric Euclidean lines are mapped
as
\begin{equation}
\mathbb{R}\rightarrow S^{1}: e^{i\pi x}=e^{\pm i\pi \rho}
\end{equation}
by winding them onto the corresponding circles of the sphere.

The winding mapping of Euclidean space onto a sphere can be
extended to any number of dimensions. Here we are focusing mostly
on the case of Euclidean space, $\mathbb{R}^{3}$, generated
as the product $RP^{2}\times \mathbb{R}$ and also on the  case of
a 3-sphere generated as $RP^{2}\times S^{1}$. In both cases we
assume the Euclidean rigidity of straight lines and the
identification of the opposite points on a circle.
Euclidean space, $\mathbb{R}^{3}$, can be mapped into a sphere,
$S^{3}$, by the winding transformation analogous to \ref{sphere}.
Indeed, for this purpose we only have
to establish a relation between the length of the radius-vector
in Euclidean space and the spherical coordinate (latitude)
measured modulo $2\pi$. The relevant transformations are
as follows:
\begin{equation}\label{sphereS3}
  \theta_{1}=\vartheta,\quad \theta_{2}=|\varphi|\mod \pi,
  \quad \phi=|\pm\pi\rho|\mod 2\pi,
\end{equation}
where the sign is determined by the quadrant of $\varphi$.

Let  $(e_{0},e_{1})$  be an orthonormal basis on a pseudo-Euclidean
plane with coordinates  $(x_{0},x_{1})$. Let the cylindrical coordinates
of $\mathbb{R}\times S^{1}$ be $(\phi,r)$. Then the simplest mapping
of this pseudo-Euclidean plane to the cylinder would be
\begin{equation}
\phi=|\pi(x_{0}+x_{1})| \mod 2\pi, r=x_{0}-x_{1}.
\end{equation}
That is, the first isotropic line is winded here around the
cylinder's cross-section (circle) and the second line is
identified with the cylinder's element.
In this way one can make a correspondence
between any non-isotropic (having a non-zero length) vector in the
plane and a point on the cylinder. For instance, if a vector $x$
having coordinates $(x_{0},x_{1})$ forms a hyperbolic angle
$\varphi$ with the $e_{0}$ or $-e_{0}$, then
\begin{equation}\label{cylinder1}
  \phi=|\pm\pi e^{-\varphi}\rho| \mod 2\pi=|\pi(x_{0}+x_{1})| \mod 2\pi.
\end{equation}
If this vector forms the hyperbolic angle $\varphi$ with the
$e_{1}$ or $-e_{1}$, then
\begin{equation}\label{cylinder2}
  r= \pm e^{\varphi}\rho=x_{0}-x_{1},
\end{equation}
where $\varphi= -\ln \left|\frac{x_{0}+x_{1}}{\rho}\right|$;
$\rho=\left|(x_{0}+x_{1})(x_{0}-x_{1})\right|^{1/2}$.

By analogy, one can build a winding map of the pseudo-Euclidean
plane into the torus, with the only difference that in the latter
case the second isotropic line is winded around the longitudinal
(toroidal) direction of the torus.

Now let  us consider a 6-dimensional pseudo-Euclidean space
$\mathbb{R}^{6}$ with the signature $(+,+,+,-,-,-)$.
In this case the analogue to the cylinder above is the product
$\mathbb{R}^{3}\times S^{3}$, in which the component $\mathbb{R}^{3}$
is Euclidean space. In order to wind the space
$\mathbb{R}^{6}$ over the cylinder $\mathbb{R}^{3}\times S^{3}$
we have to take an arbitrary pseudo-Euclidean plane in $\mathbb{R}^{6}$
passing through the (arbitrary) orthogonal lines $x_{k}$, $x_{p}$ that
belong to two Euclidean subspaces $\mathbb{R}^{3}$ of the space
$\mathbb{R}^{6}$.
Each plane $(x_{k},x_{p})$ has to be winded onto
a cylinder with the cylindrical coordinates $(\phi_{k},r_{p})$;
the indices $k,p$ correspond to the projective space $RP^{2}$.
We can take all the possible planes and wind them over the corresponding
cylinders. The mapping transformation of the  pseudo-Euclidean
space $\mathbb{R}^{6}$ into the cylinder $\mathbb{R}^{3}\times S^{3}$
is similar to the expressions \ref{cylinder1} and \ref{cylinder2}:
\begin{equation}\label{cylinder6}
  \phi_{k}=|\pm\pi e^{-\varphi}\rho| \mod 2\pi=|\pi(x_{k}+x_{p})| \mod 2\pi,
\end{equation}
\begin{equation}\label{cylinder6_1}
  r_{p}= \pm e^{\varphi}\rho=x_{k}-x_{p}.
\end{equation}
By fixing the running index $k$ and replacing it with zero
we can get the winding map of the Minkowski space $\mathbb{R}^{4}$
into the cylinder $\mathbb{R}^{3}\times S^{1}$, which is a
particular case (reduction) of \ref{cylinder6} and \ref{cylinder6_1}.
Conversely, by winding $\mathbb{R}^{3}$ over a 3-sphere,  $S^{3}$, we
can generalise the case and derive a winding map from  $\mathbb{R}^{6}$
into $S^{3}\times S^{3}$.

Let us consider the relationship between different orthonormal
bases in the pseudo-Euclidean plane, which is winded over a
cylinder. It is known that all of the orthonormal bases in a
pseudo-Euclidean are equivalent (i.e., none of them can be chosen
as privileged). However, by defining a regular field $c$ of unit
vectors on the pseudo-Euclidean plane it is, indeed, possible to
get such a privileged orthonormal basis  $(c,c_{1})$.  In turn, a
non-uniform unitary vector field $g(x)$, having a hyperbolic angle
$\varphi(x)$ with respect to the field $c$, would induce a
non-orthonormal frame $(g'(x),g'_{1}(x))$. Indeed, if we assume
that the following equalities are satisfied:
\begin{equation}\label{scale1}
  \pi=|\pm\pi e^{-\varphi}\rho(e^{\varphi}g)| \mod 2\pi=
|\pm\pi e^{-\varphi}\rho(g')| \mod 2\pi,
\end{equation}
\begin{equation}\label{scale2}
  \pm 1= \pm e^{\varphi}\rho(e^{-\varphi}g_{1})= \pm e^{\varphi}\rho(g'_{1}),
\end{equation}
we can derive a non-orthonormal frame $(g'(x),g'_{1}(x))$ by
using the following transformation of the orthonormal frame $(g(x),g_{1}(x))$:
\begin{equation}\label{riman_base}
  g'(x)=e^{\varphi}g(x), \quad g'_{1}(x)=e^{-\varphi}g_{1}(x).
\end{equation}
Then the field $g(x)$ would induce a 2-dimensional
pseudo-Riemann manifold with a metric tensor $\{g'_{ij}\}$
(where $i,j=0,1$), which is the same as the Gram matrix
corresponding to the system of vectors $(g'(x),g'_{1}(x))$.
A unitary vector field $g(x)$ defined in the Minkowski space
winded onto the cylinder $\mathbb{R}^{3}\times S^{1}$
would induce a 4-dimensional pseudo-Riemann manifold.
Indeed, take the orthonormal frame $(g,g_{1},g_{2},g_{3})$
derived by hyperbolically rotating the Minkowski space
by the angle $\varphi(x)$ in the
plane $(g(x),c)$. Then the Gram matrix
$g'_{ij}$  ($i,j=0,1,2,3$) corresponding to the set of vectors
$\{e^{\varphi}g, e^{-\varphi}g_{1}, g_{2}, g_{3}\}$ would
be related to the metric of the pseudo-Riemann manifold. Note, that, since
the determinant of the Gram matrix is unity \cite{gram79,everitt57},
the induced metric preserves the volume. That is, the differential
volume element of our manifold is equal to the corresponding volume
element of the Minkowski space.

\section{The dynamics of the model}

As we have already mentioned in Section\,1, the dynamics of the
velocity field $u(x,\tau)$ of an ideal incompressible fluid on the
surface of a cylinder $\mathbb{R}^{3}\times S^{1}$ can be
characterised by using the minimal volume principle, i.e., by
assuming that the 4-volume of the flow through an arbitrary
3-surface $\Sigma\subset\mathbb{R}^{3}$ during the time $T$ is
minimal under some initial and boundary conditions, namely:
\begin{equation}\label{cyl_principle}
  \delta\int_{0}^{T}\int_{\Sigma}dV\wedge u(x,\tau)d\tau=0,
\end{equation}
where $dV$ is the differential volume element of a 3-surface
$\Sigma$. This is also equivalent to the minimal mass carried by
the flow through the measuring surface during a finite time interval.

In a classical approximation, by using the winding projection of
the Minkowski space into a cylinder $\mathbb{R}^{3}\times S^{1}$,
we can pass from the dynamics defined on a cylinder to the statics
in the Minkowski space. Let the global time $t$ be
parameterised by the length of the flow line of the vector field $c$
in the Minkowski space corresponding to some regular vector
field on the cylinder and let the length of a single turn around
the cylinder be $h$. Let us take in the Minkowski space a set of
orthogonal to $c$ Euclidean spaces $\mathbb{R}^{3}$ in the
Minkowski space. The distance between these spaces is equal to
$hz$, where $z\in \mathbb{Z}$. The projection of this set of spaces
into the cylinder is a three-dimensional manifold, which we shall
refer to as a global measuring surface.
Then we can make a one-to-one correspondence between the dynamical
vector field $u(x,\tau)$ and the static vector field $g(x)$,
defined in the Minkowsky space.
Thus, in a classical approximation there exists a correspondence
between the minimisation of the 4-volume of the flow $u(x,\tau)$
on the cylinder and the minimisation of the 4-volume of the static
flow defined in the Minkowski space by the vector field $g(x)$,
namely:
\begin{equation}\label{Min_principle}
  \delta\int_{0}^{x_{0}} \int_{\Sigma'}dV\wedge g(x)dx_{0}=0,
\end{equation}
where the first basis vector $e_{0}$ coincides with the vector
$c$, and the 3-surfaces, $\Sigma'$, lie in the Euclidean sub-spaces
orthogonal to the vector $c$.
Let $\{(c_{i}\})=(c_{0},c_{1},c_{2},c_{3})$ be an orthonormal
basis in $\mathbb{R}^{4}$ such that $c_{0}=c$.
Let the reference frame bundle be such that each non-singular
point of $\mathbb{R}^{4}$ has a corresponding non-orthonormal
frame $(g_{i}(x))= (g_{0},g_{1},g_{2},g_{3})$, where $g_{0}=g(x)$,
$g_{1}=c_{1}$, $g_{2}=c_{2}$, $g_{3}=c_{3}$. Let us form a matrix
$\{g_{ij}\}$ of inner products $(c_{i},g_{j})$ of the basis
vectors $\{c_{i}\}$ and the frame $\{g_{i}\}$. The absolute value
of its determinant, $\det (g_{ij})$, is equal to the volume of the
parallelepiped formed by the vectors $(g_{0},g_{1},g_{2},g_{3})$.
It is also equal to the scalar product, $(g(x),c)$. On the other
hand, the equation $(g(x),c)^{2}=|\det G(x)|$ holds for the Gram
matrix, $G(x)$, which corresponds to the set of vectors
$\{g_{i}(x)\}$ \cite{Vin}. Then, according to the principle
\ref{Min_principle}, the vector field $g(x)$ satisfies the
variational equation
\begin{equation}\label{dinamoflows}
  \delta\int_{\Omega}(g(x),c)dx^{4}=\delta
  \int_{\Omega}|\det G(x)|^{\frac{1}{2}}dx^{4}=0,
\end{equation}
where $dx^{4}$ is the differential volume element of a cylindrical
4-region $\Omega$ of the Minkowski space, having the height $T$.
The cylinder's base is a 3-surface $\Sigma$ with the boundary
condition $g(x)=c$. In order to derive the differential equation
satisfying the integral variational equation \ref{dinamoflows}, we
have to find the elementary region of integration, $\Omega$. Let
$\Delta\pi$ be an infinitesimal parallelepiped spanned by the
vectors $\Delta x_{0}, \Delta x_{1}, \Delta x_{2},\Delta x_{3}$,
with $\omega$ being a tubular neighbourhood with the base
spanned by the vectors $\Delta x_{1}, \Delta x_{2},\Delta x_{3}$.
This (vector) tubular neighbourhood is filled in with the vectors
$|\Delta x_{0}|g(x)$ obtained from the flow lines of the vector
field $g(x)$ by increasing the natural parameter (the pseudo-Euclidean
length) by the amount
$|\Delta x_{0}|$.
Then the localisation expression of the
equation \ref{dinamoflows} gives \cite{moeller61}:
\begin{equation}\label{local_flow}
  \delta\int_{\Delta\pi} |\det G(x,t)|^{\frac{1}{2}}dx^{4}=\delta\vol
  \omega=0.
\end{equation}
Since the field lines of a nonholonomy vector field $g(x)$
are nonparallel even locally, any variation of such a field
(i.e, the increase or decrease of its nonholonomicity)
would result in a non-vanishing variation of the volume
$\vol \omega$.
Conversely, in the case of a holonomy field its variations do not
affect the local parallelism, so that the holonomicity of the
field $g(x)$ appears to be the necessary condition for the zero
variation of $\vol \omega$. Given a vector field $g(x)$ with an
arbitrary absolute value, the sufficient conditions for the
vanishing variation of the volume of the tubular
neighbourhood $\omega$ are the potentiality of this field and the
harmonic character of its potential. In terms of differential
forms these conditions correspond to a simple differential
equation:
\begin{equation}\label{stat}
  d\star g(x)=0,
\end{equation}
Where $d$ is the external differential; $\star$ is the Hodge star
operator;  $g(x)=d\varphi(x)$; and $\varphi(x)$ is an arbitrary
continuous and smooth function defined everywhere in the
Minkowski space, except for the singularity points (topological
features).
Substituting the unitary holonomy field $g(x)=k(x)d\varphi(x)$ in
\ref{stat}, where $k(x)=1/|d\varphi(x)|$, we shall find that
the unitary vector field $g(x)$ must satisfy the minimum
condition for the integral surfaces of the co-vector field
dual to $g(x)$. In this case the magnitude of the scalar
quantity $\varphi(x)$ will be equal to the hyperbolic angle
between the vectors $g(x)$ and $c$. We can also note that
the potential vector field $g(x)=d\varphi(x)$ represented by the
harmonic functions $\varphi(x)$ is the solution to the following
variational equation:
\begin{equation}\label{dinskalar}
  \delta\int_{0}^{T}\int_{\Sigma}\left[\left(\frac{\partial\varphi(x,t)}{\partial
t}\right)^{2} - \nabla^{2}\varphi(x,t)\right]dx^{3}dt=0,
\end{equation}
in which $\Sigma$ is a region in Euclidean space of the
``global'' observer; the function $\varphi(x,t)$ is defined in the
Minkowski space. Thus, the stationary scalar field $\varphi(x)$
induced by a topological feature in the global space is identical
to the Newtonian gravitational potential of a mass point.

We have to bear in mind that the space of a ``real'' observer
is curved, since the line for measuring time and the surface
for measuring the flux is defined by the vector field $g(x)$,
and not by the field $c$ as in the case of the global observer.
Therefore, if we wish to derive a variational equation corresponding
to the real observer, we have to define it on the pseudo-Riemann
manifold  $M$ induced in the Minkowski space by the holonomy
field $g(x)$, whose flux is measured through the surfaces orthogonal
to its flow lines and whose flow lines serve for measuring time.
The metric on $M$ is given by the Gram matrix of four tangent
vectors, one of which corresponds to the flow line $x'_{0}(\phi)$
parameterised by the angular coordinate of the cylindrical manifold,
and the three others are tangent to the coordinate lines of the
3-surface $x'_{1}(r),x'_{2}(r),x'_{3}(r)$ parameterised by the
Euclidean length. The following variational equation
holds for an arbitrary region $\Delta M$ of $M$:
\begin{equation}\label{Eienstain}
  \delta\int_{\Delta M}g^{2}(x')dV=0
\end{equation}
(under the given boundary conditions) where $dV$ is the differential
volume element of $M$. Note that the norm of the vector $g(x)$
coincides with the magnitude of the volume-element deformation
of the pseudo-Riemann manifold, which allows making the correspondence
between our functional and that of the Hilbert-Einstein action.

Returning to the global space, let us consider some properties of
the vector field $g(x)$. Let a point in the Minkowski space has a
trajectory $X(\tau)$ and velocity $\dot{X}$. Its dynamics is
determined by the variational equation:
\begin{equation}\label{point}
  \delta\int_{0}^{T}(g(x),\dot{X})d\tau=0.
\end{equation}
The varied here is the trajectory $X(\tau)$ in the Minkowski space
where the vector field $g(x)$ is defined and where the absolute
time $\tau$ plays the role of the evolution parameter. For small
time intervals the integral equation \ref{point} can be reduced to
\begin{equation}\label{linear}
  \delta(g(x),\dot{X})=0,
\end{equation}
which is satisfied by the differential equation
\begin{equation}\label{difpoint}
  \ddot{X}=g(X).
\end{equation}
Taking the orthogonal projection
$\xi(\tau)=\pr_{\mathbb{R}^{3}}X(\tau)$ of the trajectory of a
given topological feature in Euclidean space of the global
observer, as well as the projection
$\nabla\varphi(X)=\pr_{\mathbb{R}^{3}}g(X)$ of the vector field
$g(x)$ at the point $X(\tau)$ gives a simple differential equation
\begin{equation}\label{uskor}
  \ddot{\xi}(\tau)=\nabla\varphi(x),
\end{equation}
which (as in Newtonian mechanics) expresses the fact that the
acceleration of a mass point in an external gravitational field
does not depend on the mass.

\section{Some implications}

Let us consider some implications of our model for a
real observer in a classical approximation (by the real observer
we mean the reference frame of a topological feature).
First, we can note that a real observer moving uniformly along a
straight line in the Minkowski space cannot detect the ``relative
vacuum'' determined by the vector $c$ and, hence, cannot measure
the global time $t$. By measuring the velocities of topological
features (also uniformly moving along straight lines) our
observer would find that for gauging space and time one can use an
arbitrary unitary vector field $c'$ defined on the Minkowski
space. Therefore, the observer would conclude that the notion of
spacetime should be relative.
It is seen that the real observer can neither detect
the unitary vector field $g(x)$ nor its deviations from
the vector $c$.
However, it would be possible to measure the gradient of the
scalar (gravitational) field and detect the pseudo-Riemann manifold
induced by $g(x)$.

Indeed, in order to gauge time and distances in
different points of space (with different magnitudes of the scalar field)
one has to use the locally orthonormal basis
$\{g'_{i}\}$ defined on the 4-dimensional pseudo-Riemann manifold
with its metric tensor $\{g'_{ij}\}$. Thus, for the real observer,
the deformations of the pseudo-Euclidean space could be regarded
as if induced by the scalar field. Locally, the deformations
could be cancelled by properly accelerating the mass point
(topological feature), which implies that its trajectory
corresponds to a geodesics of the manifold.

We can see that the dynamics of a topological feature in our model
is identical to the dynamics of a mass point in the gravitational
field. Indeed, the scalar field around a topological feature is
spherically symmetric. At distance $r$ from the origin the metric
will be $e^{2\varphi}dt^{2}-e^{-2\varphi} dr^{2}$, which
corresponds to the metric tensor of the gravitational field of a
point mass, given $e^{2\varphi}\approx 1+2\varphi$ for small
$\varphi$. If $\varphi=H\tau$, i.e., hyperbolic angle $\varphi$
linearly depends from the evolutionary parameter $\tau$, then we
can compare the constant $H$ with the cosmological factor.

Let us now consider some quantum properties of our model. Let the
absolute value of the vector field $c$ be a continuous function
$|c(x)|$ in the Minkowski space. Then the angular velocity of the
flow will be:
\begin{equation}
\dot{\phi}(x)=\frac{d\phi(x)}{dt}= \frac{\pi}{h}|c(x)|,
\end{equation}
where the angular function $\phi(x)$ can be identified with the
phase action of the gauge potential in the observer space.
On the other hand, it is reasonable to associate the angular
velocity $X(\tau)$ of the topological  feature with the
Lagrangian of a point mass in the Minkowski
space:
\begin{equation}
\dot{\phi}(X)=\frac{d\phi(X)}{d\tau}=\frac{\pi}{h}L(x).
\end{equation}
Let us consider the random walk process of the topological
feature in the cylinder space $\mathbb{R}^{3}\times S^{1}$.
Let a probability density function $\rho(x)$ be defined on
a line, such that
$\rho(x)$,
\begin{equation}
\int_{-\infty}^{+\infty} \rho(x)dx=1.
\end{equation}
Let us calculate the expectation value for the random variable
$e^{i\pi x}$, which arises when a line is compactified into
a circe:
\begin{equation}
M(e^{i\pi x})=\int_ {-\infty}^{+\infty}\rho(e^{i\pi x})dx=\int_
{-\infty}^{+\infty} e^{i\pi x} \rho(x) dx=pe^{i\pi \alpha}.
\end{equation}
Here the quantity $pe^{i\pi \alpha}$ can be called the complex
probability amplitude. It characterises two parameters of the
random variable distribution, namely, the expectation value
itself,  $e^{i\pi \alpha}$, and the probability density, $p$, i.e.
the magnitude of the expectation value. If
$\rho(x)=\delta(\alpha)$, then $M(e^{i\pi x})=1\cdot e^{i\pi
\alpha}$. Conversely, if $\rho(x)$ is uniformly distributed along
the line then the expectation value is $M(e^{i\pi x})=0$. It
follows from these considerations that a distribution in
$\mathbb{R}^{3}$ of a complex probability amplitude is related to
random events in the cylinder space $\mathbb{R}^{3}\times S^{1}$.

In order to specify the trajectories $X(\tau)$ in the Minkowski
space with an external angular potential $\phi(x)$ we shall use
the procedure proposed by Feynman \cite{Feynman}. Let the
probabilistic behaviour of the topological feature be described as
a Markov random walk in the cylinder space $\mathbb{R}^{3}\times
S^{1}$. An elementary event in this space is a free passage. In
the Minkowski space such an event is characterised by two random
variables, duration, $\Delta \tau$, and the random path vector,
$\Delta X$, whose projection into Euclidean space of the absolute
observer is $\Delta \xi$. The ratio
 $\frac{\Delta \xi}{\Delta \tau}$
is a random velocity vector, $\dot{\xi}$.
On the other hand, the free passage of a topological feature
corresponds to an increment in the phase angle
$\Delta\phi(X)=\dot{\phi}(X) \Delta \tau$ (phase action)
in the cylinder space $\mathbb{R}^{3}\times S^{1}$.

Let the probability distribution of the phase action
has an exponential form, say, $\rho (\Delta \phi)=e^{-\Delta \phi}$
(neglecting the normalisation coefficient). Then, the corresponding
probability density for the random variable $e^{ i\Delta \phi}$
will be
\begin{equation}
\rho(e^{ i\Delta \phi})= e^{-\Delta \phi} e^{ i\Delta \phi}.
\end{equation}
Using the properties of a Markov chain \cite{meyn93}, we can
derive the probability density for an arbitrary number of random
walks:
\begin{equation}
\rho(e^{i \phi})= \prod_{0}^{T}e^{- \dot{\phi}d\tau} e^{i
\dot{\phi}d\tau}.
\end{equation}
To get the expectation value of the random variable
$e^{i \phi}$ we have to sum up over the all possible trajectories,
that is, to calculate the quantity
\begin{equation}
M(e^{i \phi})= \sum\prod_{0}^{T}e^{- \dot{\phi}d\tau} e^{i
\dot{\phi}d\tau}.
\end{equation}
It is known that any non-vanishing variation of the phase action
has a vanishing amplitude of the transitional probability and, on
the contrary, that the vanishing variation corresponds to a
non-vanishing probability amplitude \cite{erdelyi56,jones65,poston78}.
Then it is seen that the integral action corresponding to
the topological feature must be minimal. It follows that the
"probabilistic trap" of a random walk \cite{Feller} in the
cylinder space $\mathbb{R}^{3}\times S^{1}$ is determined by the
variational principle -- the same that determines the dynamics of
a mass point in classical mechanics.

 \section{Conclusions}

 In conclusion, we have made an attempt to describe the
 dynamics of spacetime (as well as of matter particles)
 in terms of the vector field defined on a cylindrical manifold
 and based on the principle of maximum mass carried by the
 field flow. The analysis of the observational implications
 of our model sheds new light on the conceptual problems
 of quantum gravity.

 Still many details of our model are left unexplored. For example,
 it would be instructive to devise the relationship between the
 vector field  $g(x)$ and the 4-potential of electromagnetic field
 $A(x)$ and to consider the local perturbations of $g(x)$
 as gravitons or/and photons. We also expect that
 the most important properties of our model would be revealed
 by extending it to the cylindrical manifold
 $\mathbb{R}^{3}\times S^{3}$. In particular, we hope that
 within such an extended version of our framework it would be
 possible to find a geometric interpretation of all known
 gauge fields.  It is also expected that studying the dynamics
 of the minimal unit vector field on a 7-sphere should be
 interesting for cosmological applications of our approach.

%%%%%%%%%%%%%%%%%%%%%%%%

\end{document}